\documentclass[11pt]{amsart}
\numberwithin{equation}{section}
\newtheorem{theorem}{Theorem}

\newtheorem{proposition}{Proposition}
\newtheorem{lemma}{Lemma}

\numberwithin{theorem}{section} \numberwithin{lemma}{section}
 \numberwithin{definition}{section}
\numberwithin{proposition}{section}
\newtheorem{remark}{Remark}[section]
\def\R{\bf R}

\def\al{\aligned}
\def\eal{\endaligned}
\def\M{{\bf M}}
\def\be{\begin{equation}}
\def\ee{\end{equation}}
\def\lab{\label}

\def\e{\epsilon}

\def\R{\bf R}
\def\M{\bf M}

\def\al{\aligned}
\numberwithin{equation}{section}

\begin{document}

\tracingpages 1
\title[ancient solutions]{\bf On ancient solutions of the heat equation }
\author{Fanghua Lin and Qi S. Zhang}
\address{Courant institute, New York University, New York, NY10012 USA  and Department of
Mathematics, University of California, Riverside, CA 92521, USA }
\date{December 2017-May 2018}

\begin{abstract}
An explicit representation formula for all positive ancient solutions of the heat equation in the Euclidean case is found. In the Riemannian case with nonnegative Ricci curvature, a similar but less explicit formula is also found. Here it is proven that any positive ancient solution is the standard Laplace transform of positive solutions of the family of elliptic operator $\Delta - s$ with $s>0$. Further relaxation of the curvature assumption is also possible. It is also shown that the linear space of ancient solutions of polynomial growth has finite dimension and these solutions are polynomials in time.
\end{abstract}
\maketitle
%\tableofcontents
\section{Introduction}

   The study of global solutions of a differential equation is a classical mathematical topic. For example, the Liouville theorem in $\mathbb{R}^n$ states that bounded harmonic functions or positive harmonic functions are constant. For evolution equations such as the heat equation, the corresponding notion of global solutions are ancient solutions i.e. solutions whose existence time is at least $(-\infty, 0)$. Understanding of ancient solutions is also useful for singularity analysis for nonlinear evolution equations.
  A Liouville type theorem for the heat equation,  different from the one for harmonic functions,  is proven in \cite{SZ:1}. In the case of Riemannian manifolds with nonnegative Ricci curvature, it  states that a nonnegative ancient solution to the heat equation whose growth rate is slower than $e^{C d(x)+ C \sqrt{|t|}}$ is a constant. Here $d(x)$ is the distance from $x$ to a reference point. Contrary to intuition, the condition is qualitatively sharp since $u=e^{x+t}$ is a nonconstant ancient solution in $\mathbb{R} \times (-\infty, 0]$.
   So a basic question arises:

{\it What are all positive ancient solutions to the heat equation?}

This is a part of a wider question of finding representation formulas for solutions of parabolic equations,  which has had a
long history, starting from Widder type uniqueness result for positive solutions.  More recently, in the interesting paper \cite{Mu:3}, p178 and also section 6, Murata also raised the open problem  to determine all positive solutions to the standard second order parabolic equations on $D \times I$ where $D$ is a noncompact domain in a Riemannian manifold and $I$ is a time interval. He treated many cases of domain $D$ such that $I$ has a left end or $I=\mathbb{R^-}$ ( section 6) and obtained  implicit representation formula. The main task seems to be to identify the parabolic Martin boundary at infinity. The main assumption is that the leading elliptic operator in the parabolic equation satisfies the intrinsic ultra-contractivity property (IU), which amounts to the operator having discrete spectrum. This left open the fundamental cases that $D$ is a typical noncompact manifold including $\mathbb{R}^n$.  By showing positive ancient solutions are completely monotone in negative time varibale, Theorem \ref{thhern} below gives a new representation formula for the most basic case of the problem for the heat equation,
i.e. when $D$ is $\mathbb{R}^n $and $I$ is the
left half line or the whole line.
 Next we extend this result to the case of manifolds with nonnegative Ricci curvature and beyond. But the result is somewhat less explicit. See Theorem \ref{thhebfm} and Remark \ref{re2.1} below. The result seems overdue, for an equation as useful and ubiquitous as the heat equation.

\begin{theorem}
\lab{thhern}
Let $u=u(x, t)$ be a nonnegative ancient solution to the heat equation in $\mathbb{R}^n \times (-\infty, 0]$. Then $u(x, -t)$ is a completely monotone function in $t$. Moreover
there exists a family of nonnegative Borel measure $\mu=\mu(\cdot, s)$ on the unit sphere $S^{n-1}$, and a Borel measure $\rho=\rho(s)$ on $[0, \infty)$ such that
\[
u(x, t) =\int^\infty_0 \int_{S^{n-1}} e^{t s + \sqrt{s} \, x \cdot \xi} d\mu(\xi, s) d\rho(s).
\]
\end{theorem}

We should point out that in 1963, Widder \cite{W:1} obtained a Cauchy type representation formula for all
nonnegative ancient solutions of the heat equation in $\mathbb{R}^n \times (-\infty, 0]$.
For example, if $u$ is a positive ancient solutions of the heat equation in $\mathbb{R}^n \times (-\infty, 0]$, then
\[
u(x, t) = \int_{\mathbb{R}^n} e^{x \cdot y+  t |y| ^2} d\mu(y)
\]for a nonnegative Borel measure.
His method is to invert his own representation formula for the Cauchy problem by the Appell transform. This method does seem to work on manifolds.

 In the papers \cite{KT:1} and \cite{Pi:1}, Koranyi-Taylor and Pinchover used an
interesting argument via the Harnack inequality to describe minimum
positive solutions for the heat equation on homogeneous manifolds, and
for parabolic equations in $\mathbb{R}^n$ with periodic coefficients. Indeed, the
Harnack inequality implies that the set of these positive solutions (after
normalizing the value at a given point in the space-time forms a noncompact, well capped convex set with compact caps in a Fr\'echet space. Then Choquet's theorem leads to a
representation in terms of an integral on the subset of extreme points
(minimal solutions). In fact, in \cite{LP:1} Section 7, Lin and Pinchover extended
this method to manifold case under a group action condition, and a class
of uniformly parabolic operators $L - d/dt$. They proved that the minimal
positive solutions is of the the form $e^{st} h(s, x)$ where h(s, x) is a
minimal positive solution of the elliptic equation.
\[
L h (s, x) - s h(s, x)=0.
\]  See also an earlier result by Murata \cite{Mu:4} in the Euclidean setting.

In contrast,  Theorem \ref{thhern} and Theorem \ref{thhebfm} below show, for manifolds with
nonnegative Ricci curvature, the following explicit formula holds:
positive ancient solutions of the heat equation is the Laplace transform,  under a Borel measure on $\mathbb{R_+}$,  of a family of positive solutions of elliptic equations. Namely, they are of the form
\[
\int e^{s t} h(s, x) d\rho(s)
\]where $h(s, x)$ is a positive solution of the elliptic equation
$
\Delta h(s, x) - s h(s, x)=0,
$ and $\rho$ is a Borel measure. For the proof, we make use of a new observation that ancient positive solutions are Bernstein's completely monotone functions for the negative time variable.

 Also the parabolic Martin boundary is identified as a family of  Martin boundaries of the elliptic equations $\Delta u- s u=0$ with $s>0$. Using the methods in the papers of Pinchover \cite{Pi:1}, Avellaneda and Lin \cite{AL:1}, and Lin and Pinchover \cite{LP:1}, one should be able to make the formula more explicit; and can generalize the main results of the present paper to the case when the Laplacian is replaced by
elliptic equations with periodic or more general coefficients. See further discussion on generalization at the end of Section 2.

    Next we turn to ancient solutions which may change sign. It is proven in \cite{SZ:1} that sublinear ancient solutions are constants. So the next step is to consider ancient solutions of polynomial growth. We will prove the finiteness of the dimension of the space of ancient solutions of polynomial growth on some Riemannian manifolds, including those with nonnegative Ricci curvature. The corresponding property for harmonic functions have been established by Colding and Minicozzi \cite{CM:1}, proving a conjecture by S.T. Yau. See also a different proof by P. Li \cite{Li:1} and other related works \cite{CM:2}, \cite{CM:3}, \cite{LW:2}, \cite{KP:1} and \cite{KP:2}.

In order to state the result, let us introduce first some notations.
We use $\M$ to denote a $n$ dimensional, complete, noncompact Riemannian manifold. The basic assumptions are:

A. volume doubling property: there exists a positive constant $d_0$ such that

\be
\lab{vdb}
|B(x, 2 r)| \le d_0 |B(x, r)|
\ee for all $x \in \M$ and $r>0$.

B. Mean value inequality for the heat equation. Let $u$ be a solution of the heat equation
$\Delta u - \partial_t u=0$ in ${\M} \times \mathbb{R}^-$. Then  for a positive constant $m_0$ and any $r>0$,
\be
\lab{mvip}
|u(x, t)|  \le \frac{m_0}{|B(x, r)| r^2} \int^t_{t-r^2} \int_{B(x, r)} |u(y, s)| dyds.
\ee It is known that both property hold if the Ricci curvature is nonnegative. See \cite{Li:1} e.g.

 We use $H^q({\M} \times \mathbb{R}^-)$ to denote the space of ancient solutions of the heat equation with growth rate at most $q \ge 1$.
i.e.
\be
\lab{upolyq}
|u(x, t)| \le c_0 (d(x, x_0) + \sqrt{|t|}+1)^q, \quad \forall (x, t) \in {\M} \times \mathbb{R}^-.
\ee
Here $x_0$ is a reference point on $\M$ and $c_0$ is any positive constant. Here and henceforth $\mathbb{R}^-=(-\infty, 0]$.

\begin{theorem}
\lab{thhepoly}
(a). Let $\M$ be a complete, n dimensional, noncompact Riemannian manifold on which assumptions $A$ and $B$ on the volume doubling property and mean value inequality for the heat equation hold.
Then, there are constants $C$ and $\eta$ such that
\[
dim (H^q({\M} \times \mathbb{R}^-)) \le C q^{\eta+1}.
\]Here $C$ depends only on the constants in the assumptions A and B and $\eta=\log_2 d_0$.

(b).  Under the same assumptions as in (a),  let $u$ be an ancient solution to the heat equation with polynomial growth of degree at most $q$, i.e. (\ref{upolyq}) holds, and let $k$ be the least integer greater than $q/2$. Then
\[
u(x, t)=u_0(x) + u_1(x) t +...+ u_{k-2}(x) t^{k-2} + u_{k-1}(x) t^{k-1},
\]with $\Delta u_{i}(x) = (i+1) u_{i+1}$, $i=0, ..., k-2$ and $u_{k-1}$ is a harmonic function.
\end{theorem}

It is well known that if $\M$ is a complete, n dimensional, noncompact Riemannian manifold with nonnegative Ricci curvature, then the assumptions in the theorem hold. In the Euclidean case, it is easy to show that harmonic functions of polynomial growth are polynomials. Likewise it is  known for long time that ancient solutions of polynomial growth are polynomials. See \cite{Ni:1} and \cite{Ei:1} e.g.
 A proof can be done by considering the spatial derivatives of the solution, which are also solutions. Then the parabolic mean value inequality can be applied to reach the conclusion. To our knowledge, the above theorem is the first such result for manifold case.

We will prove Theorem \ref{thhern} in section 2, using Li-Yau's gradient estimate and Bernstein's theorem on completely monotone functions. Theorem \ref{thhepoly} will be proven in Section 3 using an argument adapted from the case for harmonic functions in \cite{Li:1}, \cite{Ha:1} and some new input about time derivatives of solutions.

\section{Positive ancient solutions}

Proof of Theorem \ref{thhern}. At the first glance, one may use the Laplace transform on the heat equation to reduce the problem to time independent case. However there is extra boundary term which will complicate the situation. Instead we will first show that ancient positive solutions are completely monotone in time.

Let $u$ be a nonnegative ancient solution of the heat equation. Fix $(x, t) \in {\mathbb{R}^n} \times \mathbb{R}^-$. By the local Li-Yau gradient bound \cite{LY:1} applied on the parabolic cube
\[
Q_{R, T}(x, t) =\{(y, s)
\, | \, y \in B(x, R), \, s \in [t-T, t] \},
\] for $R, T>0$
there exists $C_n >0$ such that
\be
\lab{ly}
\left(\frac{1}{2} \frac{ |\nabla u|^2}{ u^2} - \frac{u_t}{u} \right) (x, t) \le C_n (\frac{1}{R^2} +\frac{1}{T}).
\ee Note that the coefficient $\frac{1}{2}$ in the first term can be replaced by any positive constant strictly less than $1$.
Hence
\[
 \frac{u_t}{u} (x, t) \ge -C_n (\frac{1}{R^2} +\frac{1}{T}).
\]Since $u$ is ancient solution, we can let $R, T$ going to infinity to deduce
\[
u_t \ge 0.
\]
Now that $u_t$ is also an ancient nonnegative solution, we can repeat the above argument to show
$u_{tt} \ge 0$ and
\[
\partial^k_t u \ge 0, \quad k=1, 2, 3, ....
\]

Fixing $x \in \R$, the one variable function of $t$
\be
\lab{fxt}
f^x(t)=u(x, -t), \quad t \in (0, \infty)
\ee is  Bernstein's completely monotone function  since
\[
(-1)^k \frac{d^k f^x(t)}{d t^k} \ge 0.
\]According to Bernstein, see Theorem 1.4 in \cite{SSV:1} e.g.,
\be
\lab{fxt=int0}
f^x(t)=  \int^\infty_0 e^{- t s} d \nu(s, x)
\ee where $\nu(\cdot, x)$ is a nonnegative Borel measure on $[0, \infty)$. Since we also have to deal with the  variable $x$, it is helpful to modify the above formula so that the measure $\nu$ can be replaced by a function. First we rewrite (\ref{fxt=int0}) as
\[
f^x(t)= a(x)  + \int^\infty_0 (e^{- t s}-1) d \nu(s, x),
\]where $a(x)=f^x(0)=u(x, 0)$. Then using Fubini theorem we compute
\[
\al
f^x(t)&=a(x)  + \int^\infty_0 (-t) \int^s_0 e^{-t \lambda} \, d\lambda d \nu(s, x)\\
&=a(x)  - \int^\infty_0  t e^{-t \lambda} \int^\infty_\lambda d \nu(s, x)  \, d\lambda \\
&= \int^\infty_0  t e^{-t \lambda} \left[a(x) - \int^\infty_\lambda d \nu(s, x)\right]  \, d\lambda
\eal
\]Define
\be
\lab{msx}
h(x, \lambda) \equiv a(x) - \int^\infty_\lambda d \nu(s, x)=\int^\lambda_0 d \nu(s, x).
\ee Then $h(x, \cdot)$ is a right continuous, non-decreasing function and
after renaming the variable $\lambda$ by $s$, we deduce
\be
\lab{fxt=int}
f^x(t)= \int^\infty_0  t e^{-t s} h(x, s) \, ds.
\ee
Fixing $x$, this shows
\be
\lab{e-ts=ux-t}
\int^\infty_0   e^{-t s} h(x, s) \, ds = \frac{u(x, -t)}{t}=\frac{u(x, -t)-u(x, 0)}{t} + \frac{u(x, 0)}{t}, \qquad t \in (0, \infty).
\ee Recall that $u=u(\cdot, \cdot)$, as a nonnegative  ancient solution is nondecreasing in time and smooth. Hence,
for $x$ in a compact set,  the function $\frac{u(x, 0)-u(x, -t)}{t}$ and derivatives are uniformly  bounded for $t\in [0, \infty)$, i.e. for each nonnegative integer $k, l$ and a compact set $D \subset \mathbb{R}^n$, there is a positive constant $C_{k, l, D}$ such that
\be
\lab{u-u/t}
\left|\partial^k_t \nabla^l \frac{u(x, 0)-u(x, -t)}{t}\right| \le C_{k, l, D}
\ee for $x \in D$ and $t \in [0, \infty)$. Since $u(x, -s)$ is completely monotone, according to Proposition 3.5 in \cite{SSV:1}, the functions (in the $t$ variable) in (\ref{e-ts=ux-t})  can be extended to the right complex plan. Therefore,
the inverse Laplace transform gives us
\be
\lab{fanlaph}
h(x, t) = \frac{1}{2 \pi i} \int^{-1 + i \infty}_{-1-i \infty} e^{s t} \frac{u(x, -s)}{s} ds.
\ee This shows that $h$ is a measurable function in $x$ and $t$.

 Observe that $ f^x(t)=u(x, -t)$ is a solution to the backward heat equation in ${\mathbb{R}^n} \times [0, \infty)$. i.e.
 \[
 \Delta f^x(t) + \partial_t f^x(t) =0.
 \]
  Let $\phi \equiv \phi(x) $ be a smooth, compactly supported test function on $\mathbb{R}^n $.
  Then
\[
\int f^x(t) \Delta \phi(x) dx + \partial_t \int f^x(t) \phi(x) dx =0.
\]This and (\ref{fxt=int0}) imply:
\[
\int \int^\infty_0   e^{- t s}  \Delta \phi(x) d\nu(x, s) -  \int \int^\infty_0  s e^{- t s}   \phi(x) d\nu(x, s)=0.
\]We mention that above integrals are convergent since $\nu(x, \cdot)$ is a nonnegative measure such that $\int^\infty_0 d\nu(x, s)  = u(x, 0)$ which is finite for each $x$. Therefore
\be
\lab{eqfornu}
 \int^\infty_0  e^{- t s}  \int ( \Delta \phi(x) -s \phi(x) ) d\nu(x, s) =0
\ee for all $t>0$.
Formally speaking this shows, by the uniqueness of the Laplace transform on measures, that in the distribution sense, we have
\be
\lab{eqh}
(\Delta -s) \nu(x, s) =0.
\ee  By \cite{Ka:1}, see also \cite{Ko:1} and \cite{CL:1}, there exists a nonnegative Borel measure $\mu=\mu(\xi, s)$ on the sphere $S^{n-1}$ such that
\[
\nu(x, s) = \int_{S^{n-1}} e^{\sqrt{s} x \cdot \xi} d\mu(\xi, s).
\] Substituting this to (\ref{fxt=int0})
we find that
\[
f^x(t)=  \int^\infty_0  e^{- t s} \int_{S^{n-1}} e^{\sqrt{s} x \cdot \xi} d\mu(\xi, s) d\rho(s).
\]Here $\rho$ is Borel measure on the positive real line. Thus, for $t<0$,
\[
u(x, t) = f^x(-t) =   \int^\infty_0  e^{ t s} \int_{S^{n-1}} e^{\sqrt{s} x \cdot \xi} d\mu(\xi, s) d\rho(s),
\] which proves the theorem. The detail of this argument is presented in the Addendum.
\qed
\medskip

Next we discuss the case where $\mathbb{R}^n$ is replaced by a noncompact Riemannian manifold with nonnegative Ricci curvature. In this general case, the Li-Yau bound (\ref{ly}) still holds. Therefore  (\ref{eqh}) is still valid and the proof is identical.
Hence
\be
\lab{uxt=h}
u(x, t) = \int^\infty_0  e^{ t s} h(s, x) ds
\ee where $h$ solves
\be
\lab{ddh-sh}
\Delta h(s, x) - s h(s, x)=0, \qquad x \in \M.
\ee So the problem of classifying ancient positive solutions is converted to classifying positive solution of the above elliptic equation. This has been dealt with, at least in the Euclidean case by M. Murata \cite{Mu:1} who generalized Martin's [Ma] result for harmonic functions to solutions of (\ref{ddh-sh}). In fact he worked on equations which includes (\ref{ddh-sh}) as a special case. See also the paper \cite{T:1} where Martin's result is generalized to the case of second order elliptic operators on manifolds, which may not be symmetric. If one does not insist on knowing the explicit form of the Martin boundary, then the following result (Proposition \ref{prellimartin} below) for  equation (\ref{ddh-sh}) follows from Martin's original method without much extra effort. First let us introduce some terminologies.

Let $\Gamma_s=\Gamma_s(x, y)$ be the minimum Green's function of the operator $\Delta -s$.
A positive solution $u$ is called minimum if the following holds: if $v$ is another solution such that $0 \le v(x) \le u(x)$ then $v(x)=c u(x)$ for a constant $c$.
Define, for a fixed $x_0 \in \M$ and points $x, y \in \M$:
\be
\lab{depxy}
\al
P_s(x, y) =
\begin{cases}
\frac{\Gamma_s(x, y)}{\Gamma_s(x_0, y)}, \quad y \neq x_0\\
0, \quad y=x_0, x \neq y,\\
1, \quad x=y=x_0.
\end{cases}
\eal
\ee

Let $\{y_j \}$ be a sequence in $\M$, which does not have an accumulation point. It is called a fundamental sequence if the function series $\{ P_s(\cdot, y_j) \}$, converges, in $C^\infty_{loc}$ topology, to a positive solution of (\ref{ddh-sh}). Two fundamental sequences $\{y_j \}$ and
$\{y'_j \}$ are called equivalent if
\[
\lim_{j \to \infty} P_s(x, y_j) = \lim_{j \to \infty} P_s(x, y'_j), \quad \forall x \in \M.
\] The symbol $\bf \Sigma$ will denote the set of equivalent classes of all fundamental sequences.  Given $w \in \Sigma$, define
\[
P_s(x, w) = \lim_{j \to \infty} P_s(x, y_j),
\]where $\{ y_j \}$ is a fundamental sequence representing $w$. For $z, z' \in \M \cup \Sigma$, a Martin distance is defined by
\[
L_s(z, z') = \int_{B(x_0, 1)} \frac{|P_s(x, z)-P_s(x, z')|}{1+|P_s(x, z)-P_s(x, z')|} dx.
\]

\begin{proposition} (c.f. Murata \cite{Mu:1} Theorem 2.3)
\lab{prellimartin}
(i) $L$ is a metric on $\M \cup \Sigma$, which is compact under $L$, and $\mathbf{\Sigma}$ is the boundary of $\M \cup \Sigma$.

(ii) $P_s(x, z)$ is continuous on $\M \times (\M \cup \Sigma)$ except when $x=z$.

(iii) Any minimal solution of (\ref{ddh-sh}) is equal to $P(x, w)$ for some $w \in \mathbf{\Sigma}$.

(iv) The set
\[
\mathbf{\Sigma_0} \equiv \{ w \in \mathbf{\Sigma} \, | \, P(\cdot, w) \, \text{is a minimal solution} \, \}
\]is a countable intersection of open sets in $\mathbf{\Sigma}$.

(v) For any positive solution $h$ of (\ref{ddh-sh}), there exists a unique nonnegative Borel measure $\mu$ on
$\mathbf{\Sigma}$ such that $\mu(\mathbf{\Sigma}-\mathbf{\Sigma_0})=0$ and
\[
h(x) =\int_{\mathbf{\Sigma_0}} P_s(x, w) d\mu(w).
\]
\end{proposition}

Based on this proposition and (\ref{uxt=h}), we immediately deduce:

\begin{theorem}
\lab{thhebfm}
Let $u=u(x, t)$ be a nonnegative ancient solution to the heat equation in ${\M}^n \times (-\infty, 0]$, where $\M$ is a complete noncompact Riemannian manifold with nonnegative Ricci curvature. Then $u(x, -t)$ is a completely monotone function in $t$. Moreover
there exists a family of nonnegative Borel measure $\mu=\mu(\cdot, s)$ on the Martin boundary
$\mathbf{\Sigma_s}$ of equation (\ref{ddh-sh}), and a Borel measure $\rho=\rho(s)$ on $[0, \infty)$, such that
\[
u(x, t) =\int^\infty_0 \int_{\mathbf{\Sigma_s}} e^{t s} P_s(x, w) \, d\mu(w, s) d\rho(s).
\]
\end{theorem}

Under some natural conditions on the manifold $\M$ in terms of the Green's function, one can prove that the Martin boundaries for
different parameters $s$ are equivalent. See p180 \cite{Mu:3}. Also it would be interesting if one can identify certain conditions on $\M$ to yield a more explicit formula for the minimal functions $P_s(x, z)$. It is well known that when the sectional curvatures are bounded between two negative constants, Anderson and Schoen \cite{AS:1} and Ancona \cite{An:1} , have worked out interesting extensions of Martin's theorem for harmonic functions. When the Ricci curvature is nonnegative, the corresponding theorem for harmonic functions is trivial since they are constants due to Liouville theorem. In instead of harmonic functions, the correct functions to study seem to be solutions of
$\Delta h- s h=0$, $s>0$. Indeed, the set of all positive entire solutions , which are normalized to be 1 at a reference point in a complete Riemannian manifold, to such an elliptic equation is convex
and compact. The compactness follows from the Moser's Harnack inequality. Thus a general representation of points (solutions) in such a compact convex in terms of an integral with a Radon measure defined on the subset of extreme points for this convex set.The latter is the statement of the classical Choquet's theorem. It is interesting however, we identified these extreme sets by its corresponding Martin boundaries. Result about these positive solutions lead to understanding of ancient solutions as we have shown.

\begin{remark}
\lab{re2.1}
The extensions of Theorem \ref{thhern} or refinement of Theorem \ref{thhebfm} to manifolds with weaker curvature condition
is also possible. For example, as indicated in \cite{Mu:4}, \cite{LP:1} and \cite{Pi:1},  suppose a uniformly restricted Harnack inequality
\be
\lab{URPH}
u(x, t-\tau) \le C_{\tau} u(x, t)
\ee holds for positive solutions of the heat equation; here $C_{\tau}$ is a positive constant depending on $\tau$ but not on $x, t$. Then minimal positive ancient solutions are of the form $e^{s t} h(s, x)$ and hence increasing in time. The integration, in the spirit of Choquet theory, of these minimal solutions is also increasing in time. This tells us that ancient positive solutions $u=u(x, t)$ are monotone nondecreasing in time. Consequently $\partial_t u$ is also a nonnegative ancient solution. Repeating this process, we find that $\partial^k_t u \ge 0$ for all $k=0, 1, 2, ...$. Now we can rerun the arguments in
Theorem \ref{thhern} or Theorem \ref{thhebfm} to conclude that $u$ is the Laplace transform of
positive solutions of the elliptic equations $\Delta h- s h=0$. It is well known from \cite{LY:1} that (\ref{URPH}) holds if the Ricci curvature is bounded from below by a negative constant and
many relaxation of this using Poincar\'e, Sobolev inequalities and volume doubling condition exist \cite{Gr:1} and \cite{Sa:1}.  As one application, with a little more work, the elliptic results \cite{AS:1} and \cite{An:1} can be extended to the parabolic case as in Theorem \ref{thhebfm} with a more explicit Martin boundary
$\Sigma_s$. As pointed out by the referee, the same techniques also apply to nonnegative solutions of parabolic equations of the form $\partial_t u
+ L u -\lambda_0 u =0$, where $L$ is a time independent elliptic operator and $\lambda_0$ is its generalized principal eigenvalue; in addition, uniqueness of nonnegative solutions of the Cauchy problem can also be derived.

\end{remark}

\section{Ancient solutions of polynomial growth}

In this section we will prove Theorem \ref{thhepoly}. First we need some notations and two lemmas.

Let $(x, t), (y, s) \in \M \times (-\infty, 0)$. One defines the parabolic distance
\[
d_p( (x, t), (y, s) ) = d(x, y) + \sqrt{|t-s|}.
\]It is easy to check that the triangle inequality holds and this is indeed a distance function. Given $(x, t), \in \M \times (-\infty, 0)$ and $r>0$, we will work with the truncated paraboloid
\[
P_r(x, t) = \{ (y, s) \, | \, d_p( (x, t), (y, s) ) \le r, \quad s \le t \}.
\]  By the doubling assumption A with doubling constant $d_0$, it is well known that for $\eta=\log_2 d_0$, and $r_2>r_1$,
\be
\lab{vdouble}
|B(x, r_2)| \le d_0 (r_2/r_1)^\eta |B(x, r_1)|
\ee

\begin{lemma}
\lab{leintu2i<}
Let $K$ be a $k$ dimensional subspace of $H^q({\bf M} \times \mathbb{R}^-)$ and $\{ u_i \}^k_{i=1}$ be any basis of $K$. Given $(x_0, t_0) \in {\M} \times \mathbb{R}^-$, $R \ge 1$, $\e \in (0, 1]$, the following inequality holds
\[
\Sigma^k_{i=1} \int_{P_R(x_0, t_0)} u^2_i dxdt \le C(d_0, m_0) \e^{-(\eta+1)}
\sup_{u \in <A, U>} \int_{P_{(1+\e)R} (x_0, t_0)} u^2 dxdt,
\]where
\[
<A, U> = \{ v = \Sigma^k_{i=1} a_i u_i \, | \, \Sigma^k_{i=1} a^2_i =1, \, a_i \in [0, 1] \}.
\]Here the constant $C(d_0, m_0)$ depends only on $n$ and the constants $d_0, m_0$ in assumptions A and B.
\end{lemma}

\proof For one fixed $(x, t) \in P_R(x_0, t_0)$, one can find $a_1, ..., a_k \in [o, 1]$ with
$\Sigma^k_{i=1} a^2_i=1$, such that
\[
\Sigma^k_{i=1} u^2_i(x, t) = v^2(x, t)
\]where
\[
v= \Sigma^k_{i=1} a_i u_i.
\]Since $d_p$ is a distance, we know, for
\be
\lab{rxt}
r(x, t) \equiv R(1+\e) - d_p((x, t), (x_0, t_0)),
\ee the following holds
\be
\lab{pxtpx0t0}
P_{r(x, t)}(x, t) \subset P_{(1+\e) R} (x_0, t_0).
\ee

Under our condition on the manifold, the $L^2$ mean value inequality holds for solutions of the heat equation on the standard parabolic cubes of size $r>0$:
\[
Q_r(x, t) = \{ (y, s) \, | \, d(x, y)<r, \,
s \in [t-r^2, t] \}.
\] i.e.
\be
\lab{mviQ}
v^2(x, t) \le \frac{C}{|Q_r(x, t)|} \int_{Q_r(x, t)} v^2(y, s) dyds.
\ee Observe that the volume of $P_r(x, t)$ and $Q_r(x, t)$ are comparable:
\[
d^{-1}_0 2^{-(\eta+2)} \le |P_r(x, t)|/|Q_r(x, t)| \le 1.
\]
Indeed,
\[
|P_r(x, t)|= \int_{B(x, r)} \int^t_{t-(r-d(x, y))^2} ds dy =\int_{B(x, r)}(r-d(x, y))^2 dy.
\]Hence
\[
\al
r^2 |B(x, r)| &\ge |P_r(x, t)| \ge \int_{B(x, r/2)} (r-d(x, y))^2 dy\\
&\ge \frac{r^2}{4} \int_{B(x, r/2)}  dy = \frac{r^2}{4} |B(x, r/2)| \ge d^{-1}_0 2^{-(\eta+2)} r^2 |B(x, r)|,
\eal
\]where the volume doubling property has been used. Therefore, (\ref{mviQ}) implies, for all $r>0$,
\be
\lab{mviP}
v^2(x, t) \le \frac{C}{|P_r(x, t)|} \int_{P_r(x, t)} v^2(y, s) dyds.
\ee Here the constant $C$ may have changed.

The above calculation also implies a volume comparison result for the paraboloid $P_r(x, t)$. For $r_2>r_1>0$, we have
\be
\al
\lab{pr2pr1}
\frac{|P_{r_2}(x, t)|}{r^{\eta+2}_2} \le \frac{|B(x, r_2)|}{r^\eta_2} \le C \frac{|B(x, r_1)|}{r^\eta_1}
\le C \frac{|P_{r_1}(x, t)|}{r^{\eta+2}_1}.
\eal
\ee In (\ref{pr2pr1}), taking $r_2=2 R$ and
\[
r_1=r(x, t) = R(1+\e) - (d(x, x_0) + \sqrt{|t-t_0|}),
\]we deduce
\[
\al
|P_r(x, t)| &\ge C^{-1} |P_{2 R}(x, t)| \left[ \frac{R(1+\e) - (d(x, x_0) + \sqrt{|t-t_0|})}{2R} \right]^{\eta+2}\\
&\ge C^{-1} [R(1+\e) - (d(x, x_0) + \sqrt{|t-t_0|})]^{\eta+2} \frac{|P_R(x_0, t_0)|}{R^{\eta+2}}.
\eal
\]Here we have used the volume doubling property to shift the vertex of the paraboloid.
Substituting the above inequality into (\ref{mviP}), we find that
\[
\al
v^2(x, t) &\le \frac{C}{[[R(1+\e) - (d(x, x_0) + \sqrt{|t-t_0|})]^{\eta+2}} \int_{P_r(x, t)} v^2(y, s) dyds \frac{R^\eta}{|B(x_0, R)|}\\
&\le \frac{C}{[[R(1+\e) - (d(x, x_0) + \sqrt{|t-t_0|})]^{\eta+2}} \int_{P_{R(1+\e)}(x_0, t_0)} v^2(y, s) dyds \frac{R^\eta}{|B(x_0, R)|}.
\eal
\]Here we just used the relation (\ref{pxtpx0t0}).

Next we integrate the last inequality on $P_R(x_0, t_0)$ to deduce
\be
\lab{Su2i<}
\al
\int_{P_R(x_0, t_0)} \Sigma^k_{i=1} u^2_i(x, t)dxdt &\le C \int_{P_R(x_0, t_0)} \frac{dxdt}{[R(1+\e) - (d(x, x_0) + \sqrt{|t-t_0|})]^{\eta+2}} \\
&\qquad \qquad \times \sup_{v \in <A, U>} \int_{P_{R(1+\e)}(x_0, t_0)} v^2 dxdt \,
\frac{R^\eta}{|B(x_0, R)|}.
\eal
\ee To finish the proof of the lemma, we compute the first integral on the right hand side of the above inequality.
\[
\al
I &\equiv \int_{P_R(x_0, t_0)} \frac{dxdt}{[R(1+\e) - (d(x, x_0) + \sqrt{|t-t_0|})]^{\eta+2}}\\
&=\int_{B(x_0, R)} \int^{(R-d(x, x_0))^2}_0 \frac{dt}{[R(1+\e) - (d(x, x_0) + \sqrt{t})]^{\eta+2}} dx.
\eal
\]Writing $A(x)=R-d(x, x_0)$ and making the change of variables
$
s = \sqrt{t}, \quad dt= 2s ds,
$ we deduce
\[
\al
I &= \int_{B(x_0, R)} \int^{A(x)}_0 \frac{2 s ds}{[R \e + A(x) - s]^{\eta+2}} dx\\
&\le \int_{B(x_0, R)} 2 A(x) \int^{A(x)}_0 \frac{ ds}{[R \e + A(x) - s)]^{\eta+2}} dx\\
&\le \frac{2}{n+1} \int_{B(x_0, R)} 2 A(x) \frac{1}{[R \e]^{\eta+1}} dx.
\eal
\]Hence
\[
I \le C \e^{-(\eta+1)} \frac{|B(x_0, R)|}{R^{\eta}}.
\]Substituting this to (\ref{Su2i<}), we get
\[
\Sigma^k_{i=1} \int_{P_R(x_0, t_0)} u^2_i(x, t) \le C \e^{-(\eta+1)} \sup_{v \in <A, U>} \int_{P_{R(1+\e)}(x_0, t_0)} v^2 dxdt. \qed
\]

The next Lemma is a parabolic version of Lemma 3.4 in \cite{Ha:1}. We present a proof for completeness.

\begin{lemma}
\lab{le2.2} Let $K$ be a finite dimensional subspace of $H^q({\M} \times \mathbb{R}^-)$. There
exists a constant $R_0=R_0(K)$ such that, for all $R \ge R_0$, all $(x_0, t_0) \in {\M} \times \mathbb{R}^-$,
\be
\lab{uv}
<u, v> \equiv \int_{P_R(x_0, t_0)} u \, v dx dt
\ee is an inner product on $K$.
\proof
\end{lemma}

Let $\{ u_i \}^k_{i=1}$ be any basis of $K$. Suppose the lemma is not true. Then there exists a sequence $R_j \to \infty$ and $v_j = \Sigma^k_{i=1} a^i_j u_i \neq 0$  such that $\int_{P_{R_j}(x_0, t_0)} v^2_j=0$. Let us assume without loss of generality that $\Sigma^k_{i=1} (a^i_j)^2 =1$ for each $j$. Since the unit sphere $S^{k-1} \subset \R^k$ is compact, we can find a subsequence of $\{ a^i_{j} \}^\infty_{j=1}$, identically denoted, such that
\[
a^i_j \to b^i, \quad j \to \infty.
\]Write $v= \Sigma^k_{i=1} b^i u_i$. Note that $\Sigma^k_{i=1} (b^i)^2 =1$ and
$\{ u_i \}^k_{i=1}$ is a basis. Hence $v$ is not identically $0$.  However, for any fixed $R>0$, we have
\[
0=\lim_{j \to \infty} \int_{P_R(x_0, t_0)} v^2_j dxdt = \int_{P_R(x_0, t_0)} v^2 dxdt.
\] Hence $v \equiv 0$ on ${\M} \times \mathbb{R}^-$. This is a contraction. \qed

Using the techniques in \cite{Ln:1} and \cite{Po:1}, together with a new argument, one can strengthen the Lemma to show that (\ref{uv}) is an inner product for any $R>0$. To avoid interrupting the flow of the proof of the main results, we will do it in Proposition \ref{neijiR} at the end of the paper.
\medskip

Now we can prove the Theorem \ref{thhepoly}, part (a).

Let $K$ be a finite dimensional subspace of $H^q({\M} \times \mathbb{R}^-)$ and $\{ u_i \}^k_{i=1}$ be an orthonormal basis of $K$ with respect to the inner product
\[
A_{\beta R} (u, v) = \int_{P_{\beta R}(x_0, t_0)} u v dxdt.
\] Here $(x_0, t_0) \in {\M} \times \mathbb{R}^-$ and $\beta>1$ is to be chosen later. We claim that given any $\delta>0$, $R_0 \ge 1$, there exists $R>R_0$ such that
\be
\lab{intu2i>}
\Sigma^k_{i=1} \int_{P_R(x_0, t_0)} u^2_i dxdt \ge k \beta^{- (2 q + \eta+2 + \delta)}.
\ee

In order to prove the claim, let us introduce some notations.  Denote by $tr_{R'} A_R$ (respectively $det_{R'} A_R$), the trace (respectively determinant) of the inner product $A_R$ with respect to $A_{R'}$.  Thus
$tr_{R'} A_R$ (respectively $det_{R'} A_R$) is the trace (determinant) of the matrix
\[
< A_R(v_i, v_j)>_{i, j=1, ...k}
\]where $\{ v_i \}^k_{i=1}$ is an orthonormal basis of $K$ with respect to the inner product $A_{R'}$.

Suppose the claim is not true. Then, given $\delta>0$ and $\beta>1$,  there exists $(x_0, t_0) \in {\M} \times \mathbb{R}^-$, $R_0 \ge 1$ such that the following holds for all $R \ge R_0$:
\[
tr_{\beta R} A_R = \Sigma^k_{i=1} \int_{P_R(x_0, t_0)} u^2_i dxdt < k \beta^{- (2 q + \eta+2 + \delta)},
\]where $\{ u_i \}^k_{i=1}$ is an orthonormal basis of $K$ with respect to the inner product $A_{\beta R}$.
Then
\[
\left( det_{\beta R} A_R \right)^{1/k} \le \frac{ tr_{\beta R} A_R}{k}  <\beta^{- (2 q + \eta+2 + \delta)}.
\]This shows, since $det_{\beta R} A_R= (det_R A_{\beta R})^{-1}$, that
\[
det_R A_{\beta R} > \beta^{k (2 q + \eta+2 + \delta)}.
\]Iterating this for $R, \beta R, ..., \beta^j R$, we deduce
\[
det_R A_{\beta^j R} > \beta^{j k (2 q + \eta+2 + \delta)}.
\]By the growth assumption on $u_i$ and the volume bound on the paraboloid, it is easy to see that
\[
det_R A_{\beta^j R} \le k! C^k (\beta^j R)^{(2 q + \eta+2)} |B(x_0, 1)|.
\]The last two inequalities contradict each other when $j$ is sufficiently large. This shows the claim ((\ref{intu2i>})) is true.

 In Lemma \ref{leintu2i<} we take $\beta$ to be $1+\e$; then  from (\ref{intu2i>}) and we find that
\[
k \beta^{- (2 q +\eta+2 + \delta)} \le \Sigma^k_{i=1} \int_{P_R(x_0, t_0)} u^2_i dxdt \le C(d_0, m_0) \e^{-(\eta+1)}
\sup_{u \in <A, U>} \int_{P_{(1+\e)R} (x_0, t_0)} u^2 dxdt,
\]which implies
\[
k \le C (1+\e)^{ (2 q + \eta+2 + \delta)}  \e^{-(\eta+1)}.
\]Taking $\e= 1/q$ and letting $\delta \to 0$, we deduce
\[
k \le C q^{\eta+1}.
\] This proves part (a) of the theorem.

Next we prove part (b). Since the mean value inequality holds (Assumption B), from the work of Grigoryan \cite{Gr:1}, (see also Li-Wang \cite{LW:1}) and the volume doubling inequality, we know that the heat kernel has a global Gaussian upper bound
\be
\lab{gub}
G(x, t, y) \le  \frac{c_1}{|B(x, \sqrt{t})|}
e^{-c_2 \frac{d^2(x, y)}{t}}
\ee for all $x, y \in \M$ and $t>0$.

 Using this upper bound, we observe that an ancient solution of polynomial growth is also an eternal solution, namely the time of existence is $(-\infty, \infty)$. One can just use $u(x, 0)$ as the initial value to solve the heat equation. Since $u(x, 0)$ has a polynomial bound  and the heat kernel has exponential decay on space (\ref{gub}), the solution exists on $[0, \infty)$. It easy to see, via integration by parts with test functions,  that this forward solution and the ancient solution together form a smooth eternal solution by checking that at $t=0$ the solution is smooth. We also denote this eternal solution by $u=u(x, t)$.
The rest of the proof is divided into 3 steps.

{\it Step 1.} We show that the forward solution also has the polynomial bound, i.e, for a positive constant $c_0$,
\be
\lab{ut+jie}
|u(x, t)| \le c_0 (d(x, x_0) + \sqrt{t}+1)^q,  \qquad \forall (x, t) \in {\M} \times [0, \infty).
\ee

Since $|u(x, 0)| \le c_0 (d(x, x_0) +1)^q$ by assumption, for $t \ge 0$ and $x \in \M$, we have
\[
\al
|u(x, t)| &= \left| \int G(x, t, y) u(y, 0) dy \right| \le c_0 \int G(x, t, y) (d(y, x_0) +1)^q dy\\
&\le c_0 \int G(x, t, y) (d(x, x_0) + d(x, y) +1)^q dy\\
&\le c_0 2^q (d(x, x_0) +1)^q + c_0 2^q \int G(x, t, y)  d(x, y)^q dy.
\eal
\]Here $G$ is the heat kernel and in the last step the inequality $\int G(x, t, y) dy \le 1$ is used.

By (\ref{gub}),
\[
\al
|u(x, t)| &\le c_0  2^q (d(x, x_0) +1)^q + c_0 c_1 2^q \int   \frac{1}{|B(x, \sqrt{t})|}
e^{-c_2 \frac{d^2(x, y)}{t}} d(x, y)^q dy\\
&=c_0  2^q (d(x, x_0) +1)^q + c_0 c_1 2^q t^{q/2} \int   \frac{1}{|B(x, \sqrt{t})|}
e^{-c_2 \frac{d^2(x, y)}{2t}} e^{-c_2 \frac{d^2(x, y)}{2t}} \frac{d(x, y)^q}{t^{q/2}} dy\\
&\le c_0  2^q (d(x, x_0) +1)^q + c_3 t^{q/2} \Sigma^\infty_{-\infty}
\int_{2^i \sqrt{t} \le d(x, y) \le 2^{i+1} \sqrt{t}}   \frac{1}{|B(x, 2^{-i-1} d(x, y))|}
e^{-c_2 \frac{d^2(x, y)}{2t}} dy\\
&\le c_0  2^q (d(x, x_0) +1)^q + c_4 t^{q/2}.
\eal
\]In the above have just used the volume doubling condition. This proves (\ref{ut+jie}), where the constant $c_0$ may have changed.
\medskip

{\it Step 2.} We show that the $k-th$ time derivative of $u$ is $0$ for $k >q/2$.
\medskip

Fix a point $(x_1, t_1) \in \M \times \mathbb{R}$ and $R>0$, let $Q^0_R$ be the full parabolic cube
\[
\{ (x, t) \, | \, d(x, x_1)<R, \quad |t-t_1| \le R^2 \}.
\]Denote by $\psi$ a standard smooth cut off function supported in $Q^0_{2R}$ such that $\psi=1$ in $Q^0_{3R/2}$ and $|\nabla \phi|^2 + |\partial_t \psi| \le C/R^2$. Since $u$ is a smooth solution to the heat equation, we compute
\[
\al
\int_{Q^0_{2R}} &(\Delta u)^2 \psi^2 dxdt = \int_{Q^0_{2R}} u_t \Delta u \psi^2 dxdt\\
&=-\int_{Q^0_{2R}} ((\nabla u)_t \nabla u) \, \psi^2 dxdt - \int_{Q^0_{2R}} u_t \nabla u \nabla \psi^2 dxdt\\
&= - \frac{1}{2} \int_{Q^0_{2R}} (|\nabla u |^2)_t \, \psi^2 dxdt - 2 \int_{Q^0_{2R}} u_t \psi \nabla u \nabla \psi dxdt\\
&\le \frac{1}{2} \int_{Q^0_{2R}} |\nabla u |^2 \, (\psi^2)_t dxdt +\frac{1}{2} \int_{Q^0_{2R}} (u_t)^2 \psi^2 dxdt + 2 \int_{Q^0_{2R}} |\nabla u|^2 |\nabla \psi|^2 dxdt.
\eal
\]This and the standard Cacciopoli inequality (energy estimate) show that
\be
\lab{ddu2r-4}
\int_{Q^0_{R}} (\Delta u)^2  dxdt \le \frac{C_0}{R^4} \int_{Q^0_{2R}} u^2  dxdt.
\ee Here $C_0$ is a universal constant.

Since $u_t$ is also an eternal solution, we can replace $u$ in (\ref{ddu2r-4}) by $u_t$ to deduce
\[
\al
\int_{Q^0_{R}} (\partial^2_t u)^2  dxdt &= \int_{Q^0_{R}} (\Delta  u_t)^2  dxdt\\
&\le \frac{C_0}{R^4} \int_{Q^0_{2R}} (u_t)^2  dxdt \le \frac{C^2_0}{R^8} \int_{Q^0_{4R}} u^2  dxdt.
\eal
\]By induction, we deduce
\be
\lab{utk<}
\int_{Q^0_{R}} (\partial^k_t u)^2  dxdt \le \frac{C^k_0}{R^{4 k}} \int_{Q^0_{2^k R}} u^2  dxdt.
\ee

Applying the mean value inequality (A) on $\partial^k_t u$, which is also a solution to the heat equation, we find, using (\ref{utk<}), that
\[
\al
|\partial^k_t &u(x_1, t_1)|^2 \le \frac{2 m_0}{|Q^0_R|} \int_{Q^0_{R}} (\partial^k_t u)^2  dxdt
\le \frac{2 m_0 C^k_0}{R^{4 k} |Q^0_R|} \int_{Q^0_{2^k R}} u^2 dxdt\\
&\le \frac{2 m_0 C^k_0}{R^{4 k} |Q^0_R|} \int_{Q^0_{2^k R}} c^2_0 (d(x, x_0) + \sqrt{|t|} + 1)^{2q} dxdt\\
&\le \frac{2 m_0 C^k_0 c^2_0}{R^{4 k} |Q^0_R|} \int_{Q^0_{2^k R}}  (d(x, x_1) + \sqrt{|t-t_1|}+ d(x_1, x_0) + \sqrt{|t_1|} + 1)^{2q} dxdt\\
&\le \frac{2 m_0 C^k_0 c^2_0}{R^{4 k} |Q^0_R|} |Q^0_{2^k R}|  (  2^{k+1} R+ d(x_1, x_0) + \sqrt{|t_1|} + 1)^{2q}.
\eal
\] Using the volume doubling property, we arrive at
\[
|\partial^k_t u(x_1, t_1)|^2 \le \frac{2 m_0 C^k_0 c^2_0 2^{k (n+2)}}{R^{4 k}}  (  2^{k+1} R+ d(x_1, x_0) + \sqrt{|t_1|} + 1)^{2q}.
\]Letting $R \to \infty$ and using $k>q/2$, we see that $\partial^k_t u(x_1, t_1)=0$. This completes step 2.
\medskip

{\it Step 3.} From Step 2, we know that
\[
u(x, t) = u_0(x) + u_1(x) t +...+ u_{k-1}(x) t^{k-1}.
\]Substituting this to the heat equation, we deduce
\[
\Delta u_0(x) + \Delta u_1(x) t +...+ \Delta u_{k-1}(x) t^{k-1} = u_1(x) + u_2(x) 2t + ... +
u_{k-1}(x) (k-1) t^{k-2}.
\]This implies
\[
\Delta u_{k-1}(x) =0, \quad \Delta u_{k-2}(x) = (k-1) u_{k-1}, \quad ..., \quad \Delta u_0(x)=u_1(x).
\] So finally we deduce
\be
\lab{u=u0+uk}
u(x, t)=u_0(x) + u_1(x) t +...+ u_{k-2}(x) t^{k-2} + u_{k-1}(x) t^{k-1},
\ee with $\Delta u_{i}(x) = (i+1) u_{i+1}$, $i=0, ..., k-2$. This completes the proof of the theorem. \qed
\medskip

\noindent{\it Remark.} {\it If $q/2$ is an integer, then $k-1=q/2$. Fixing $x \in \M$ and $t \neq 0$ and dividing both sides of (\ref{u=u0+uk}) by $|t|^{q/2}$, we get
\[
\frac{u(x, t)}{|t|^{q/2}} = \frac{u_0(x)}{|t|^{q/2}} +...+ \frac{(-1)^{k-2} u_{k-2}(x)} {|t|} + u_{k-1}(x).
\]Letting $t \to -\infty$, from the assumed bound on $u$, we see that $u_{k-1}$ is bounded.
In case $\M$ has nonnegative Ricci curvature, Yau's Liouville theorem implies $u_{k-1}$ is a constant.}

Finally, we show that some techniques from the proof of the theorems can be used to prove certain backward uniqueness result for ancient solutions and improve Lemma \ref{le2.2}.

\begin{proposition}
\lab{neijiR}
Let $\M$ be a complete, n dimensional, noncompact Riemannian manifold on which assumptions $A$ and $B$ on the volume doubling property and mean value inequality for the heat equation hold. Then the following conclusions are true.

(a). Let $u$ be an ancient solutions of the heat equation with growth rate at most $q \ge 1$.
i.e. (\ref{upolyq}) holds. Suppose $u(x, 0)=0$. Then $u \equiv 0$.

(b).  Let $K$ be a finite dimensional subspace of $H^q({\M} \times \mathbb{R}^-)$.  For all $R >0$, all $(x_0, t_0) \in {\M} \times \mathbb{R}^-$,
\be
\lab{uv2}
<u, v> = \int_{P_R(x_0, t_0)} u \, v dx dt
\ee is an inner product on $K$.
\end{proposition}

We comment that this result does not follow from the unique continuation result for bounded solutions in \cite{Po:1} p530 Remark. One reason is that we do not have the Hamilton's matrix Harnack inequality for the heat equation, which requires nonnegative sectional curvature and parallel Ricci curvature. The other is that our solution can be unbounded.

\proof (a). The proof follows from that of Theorem \ref{thhepoly} part (b). Indeed, since $u(x, 0)=0$, we can extend $u$ to be an eternal solution by assigning $0$ value for positive time.
Then the bound (\ref{upolyq}) holds for all space time automatically, without using heat kernel bound.  Following Step 2 in the proof of Theorem \ref{thhepoly} part (b), which uses only the mean value and volume doubling properties, we find that
\[
u(x, t) = u_0(x) + u_1(x) t +...+ u_{k-1}(x) t^{k-1},
\]where $\Delta u_i = (i+1) u_{i+1}$, $i=0, 1, 2, ... k-1$. The assumption that $u(x, 0)=0$ implies $u_0=0$. Hence $u_i=0$ for all $i$ and $u \equiv 0$. This proves part (a) of the proposition.

(b).  We prove by contradiction. Suppose for some $R>0$, (\ref{uv2}) is not an inner product.
Then there is an ancient solution $u$ of polynomial growth such that
\[
\int_{P_R(x_0, t_0)} u^2 dxdt =0,
\]but $u \neq 0$. By \cite{Ln:1}, one knows that $u(\cdot, t)=0$ for all $t \in [t_0-R^2_0, t_0]$. We remark the result and method in \cite{Ln:1} is a local one, which works for the manifold case. Now part (a) of the proposition implies $u \equiv 0$, which is a contradiction.
\qed

\section { Addendum by Thomas Swayze}

{\it Department of mathematical Sciences,  Carnegie Mellon University, Pittsburgh, PA 15213,
email: tes@andrew.cmu.edu}

\medskip

To complete the proof of Theorem \ref{thhern}, we start from equation (\ref{eqfornu})
 i.e.
\be
\lab{eqfornu2}
 \int^\infty_0  e^{- t s}  \int ( \Delta \phi(x) -s \phi(x) ) d\nu(x, s) =0
\ee for all $t>0$.
Recall from (\ref{msx}) that
\be
\lab{hnu_x}
h(x,t)=\int^t_0 d \nu(x, s) \equiv \int^t_0 d \nu_x(s) \equiv \nu_x([0, t]) \le \nu_{x}([0,\infty))=u(x,0).
\ee
Hence $h$ is locally bounded. Since $h(x,t)$ is a measurable
function of $x$ for each $t$, we may treat  $\nu$ as
a measure in $\mathbb{R}^{n}\times[0,\infty)$ in the following manner: for $E\subseteq\mathbb{R}^{n}$
and $t\ge0$:
\[
\nu(E\times[0,t])=\int_{E}h(x,t)dx.
\]
Note that $\nu(K\times[0,\infty))=\int_{K}u(x,0)dx<\infty$ for all
compact $K\subset\mathbb{R}^{n}$, so $\nu$ is locally finite.

The left hand side of (\ref{eqfornu2}) is the Laplace transform of the signed measure
$\eta_{\phi}$ on $[0,\infty)$ given by:
\[
d\eta_{\phi}(s)=\int_{\mathbb{R}^{n}}\left(\triangle\phi(x)-s\phi(x)\right)d\nu(x,s).
\]
So $\eta_{\phi}$ is identically zero, and thus:
\[
0=\eta_{\phi}([0,t])=\int_{\mathbb{R}^{n}}\left[\triangle\phi(x)\int_{0}^{t}d\nu_{x}(s)-
\phi(x)\int_{0}^{t}sd\nu_{x}(s)\right]dx.
\]Here $\nu_x$ is defined in (\ref{hnu_x}).
Integrating by parts yields:
\[
\int_{\mathbb{R}^{n}}\triangle\phi(x)h(x,t)dx=\int_{\mathbb{R}^{n}}\phi(x)\left[th(x,t)-\int_{0}^{t}h(x,s)ds\right]dx
\]
As this holds for all $C_{c}^{\infty}$ function $\phi$, it follows
that $h(\cdot,t)$ is smooth for all $t$ with:
\[
\triangle h(x,t)=th(x,t)-\int_{0}^{t}h(x,s)ds.
\]

Now let $I=(a,b]$. Then for all $x\in\mathbb{R}^{n}$, $\nu_{x}(I)=h(x,b)-h(x,a)$.
So $\nu_{x}(I)$ is a smooth function of $x$, and:
\begin{align*}
\triangle\nu_{x}(I) & =bh(x,b)-ah(x,a)-\int_{a}^{b}h(x,s)ds\\
 & =b(h(x,b)-h(x,a))-\int_{a}^{b}h(x,s)-h(x,a)ds\\
 & \le b\nu_{x}(I).
\end{align*}
A similar argument shows that $\triangle\nu_{x}(I)\ge a\nu_{x}(I)$.
In particular, $d(x)=\frac{\triangle\nu_{x}(I)}{\nu_{x}(I)}$ satisfies
$|d(x)|\le b$ on $\mathbb{R}^{n}$. So by Harnack's inequality on
the elliptic operator $L=\triangle-d$, there is a constant $C(|x|,b)$
so that:
\[
\nu_{x}(I)\le C(|x|,b)\nu_{0}(I).
\]
From this, it follows that $\nu_{x}(E)\le C(|x|,b)\nu_{0}(E)$ for
all $E\subseteq[0,b]$. As this holds for all $b$, $\nu_{x}\ll\nu_{0}$.
Furthermore, $\frac{d\nu_{x}}{d\nu_{0}}(s)\le C(|x|,b)$ whenever
$s\le b$, so $\frac{d\nu_{x}}{d\nu_{0}}(s)$ is a locally bounded
function of $s$ and $x$.

Returning to (\ref{eqfornu2}):
\begin{align*}
0 & =\int_{\mathbb{R}^{n}\times[0,\infty)}e^{-st}\left[\triangle\phi(x)-s\phi(x)\right]d\nu(x,s)\\
 & =\int_{0}^{\infty}e^{-st}\left(\int_{\mathbb{R}^{n}}\frac{d\nu_{x}}{d\nu_{0}}(s)
 \left[\triangle\phi(x)-s\phi(x)\right]dx\right)d\nu_{0}(s).
\end{align*}
Again, the last line is the Laplace transform of a signed measure,
which we now know must be identically zero. So for $\nu_{0}$-almost
every $s\in[0,\infty)$, we have that:
\[
\int_{\mathbb{R}^{n}}\frac{d\nu_{x}}{d\nu_{0}}(s)\left[\triangle\phi(x)-s\phi(x)\right]dx=0.
\]
This holds for all $\phi\in C_{c}^{\infty}$, so $\frac{d\nu_{x}}{d\nu_{0}}(s)$
is smooth and satisfies:
\[
\triangle\frac{d\nu_{x}}{d\nu_{0}}(s)=s\frac{d\nu_{x}}{d\nu_{0}}(s).
\]
Thus, as pointed out in Section 2, by \cite{Ka:1}, see also \cite{Ko:1} and \cite{CL:1}, there is a Borel
measure $\mu_{s}$ on $\mathbb{S}^{n-1}$ so that:
\[
\frac{d\nu_{x}}{d\nu_{0}}(s)=\int_{\mathbb{S}^{n-1}}e^{\sqrt{s}x\cdot\xi}d\mu_{s}(\xi).
\]
Plugging this in to the  formula (\ref{fxt=int0}):
\[
u(x,-t)=\int_{0}^{\infty}e^{-st}d\nu(s, x) =\int_{0}^{\infty}e^{-st}d\nu_{x}(s),
\]
we see that:
\[
u(x,t)=\int_{0}^{\infty}\int_{\mathbb{S}^{n-1}}e^{st+\sqrt{s}x\cdot\xi}d\mu_{s}(\xi)d\nu_{0}(s).
\]This completes the proof of Theorem \ref{thhern} after renaming the measure $\nu_0$ as $\rho$.

We remark that alternatively, one can show that the gradient of the solution $u=u(x, t)$ converges to $0$ as $t \to \infty$. Then using this fact on the inverse Laplace transform
(\ref{fanlaph}) one can prove that $h(\cdot, t)$ is a smooth function in $x$ variable and proceed to reach the same conclusion.

\medskip
{\bf Acknowledgment.} We wish to thank Professors Bobo Hua and Philippe Souplet for helpful discussions. Part of the paper was written when both authors was visiting the Shanghai Center of Mathematical Sciences and the School of Mathematics at the Fudan University. We are grateful to Professor Lei Zhen and the Center for the invitation and for their warm hospitality. We are also indebted to Professor Y. Pinchover who informed us some of his and other's related works.
Thanks also go to the anonymous referee for checking the paper carefully and suggesting a number of references and possible applications. We are also indebted to Mr. Thomas Swayze and Professor Bob Pego for pointing out an error in an earlier version of Theorem \ref{thhern} and for presenting a way to fix it in the addendum.

F.H.L acknowledges the support of NSF grant DMS-1501000 and Q.S.Z is thankful of the Simons Foundation for its support.

\bigskip

\noindent e-mails:
linf@cims.nyu.edu and qizhang@math.ucr.edu

\enddocument